\documentclass[12pt]{article}
\usepackage{amsmath}
\usepackage{amsfonts}
\usepackage{amssymb}
\usepackage{amsthm}
\usepackage{hyperref}

\usepackage{graphicx}

\def\bee{\begin{equation}}
\def\eee{\end{equation}}

\begin{document}

\thispagestyle{empty}
\centerline{}
\bigskip
\bigskip
\bigskip
\bigskip
\bigskip
\bigskip
\centerline{\Large\bf New lower bounds for two color and multicolor}
\centerline{\Large\bf Ramsey numbers}

\bigskip
\bigskip

\begin{center}

{\large \sl Robert Gerbicz}\\*[5mm]
Eotvos Lorand University\\
Email: robert.gerbicz@gmail.com\\
\bigskip
May 6, 2010
\end{center}

\bigskip
\bigskip
\bigskip

\begin{center}
{\bf Abstract}\\
\bigskip
\begin{minipage}{12.8cm}
Using cyclic graphs I give new lower bounds for two color and multicolor\\
 Ramsey numbers: $R(4,16)\ge 164$, $R(5,11)\ge 171$, $R(5,12)\ge 191$, $R(5,13)\ge 213$, $R(5,14)\ge 239$, $R(3,3,9)\ge 118$, $R(3,3,10)\ge 142$ and $R(3,3,11)\ge 158$. Improving the previous best known bounds.
\end{minipage}
\end{center}

\bigskip

\section{New bounds}

By definition: $R(k,l)$ is the smallest number such that for any complete graph on $R(k,l)$ vertices, whose edges are colored by two colors (red and blue), there exists either a red $K_k$, or a blue $K_l$. So to prove that $R(k,l)\ge n+1$ it is enough to make a good edge coloring on $n$ vertices. I'm using cyclic graphs. On $n$ vertices (number these points by $1,2,\cdots,n$) the symmetric circle coloring means that the edge of $(i,j)$ is red if $d\in S$ set, where $d=$min$(|i-j|,n-|i-j|)$, otherwise the color of the edge is blue. Similarly for multicolorings we give sets and the color of an edge depends only on the value of $d$. To read a survey about Ramsey numbers, see $\bf {[1]}$. 
\\
\\
To prove that $R(4,16)\ge 164$, let $n=163$, and the set:\\
$S=\{3,17,24,25,27,37,44,45,50,53,54,55,57,63,64,65,73,78,79,80\}.$
\\
\\
To prove that $R(5,11)\ge 171$, let $n=170$, and the set:\\
$S=\{4,5,6,7,12,17,18,19,25,26,27,28,33,36,41,42,43,44,49,53,54,\\
55,56,57,58,59,60,65,69,73,77,81,85\}.$
\\
\\
To prove that $R(5,12)\ge 191$, let $n=190$, and the set:\\
$S=\{1,2,3,5,8,11,13,18,20,22,23,27,28,33,34,37,38,41,42,43,47,\\
48,49,53,54,55,58,59,62,65,71,73,74,81,83,93,95\}.$
\\
\\
To prove that $R(5,13)\ge 213$, let $n=212$, and the set:\\
$S=\{2,4,5,13,15,16,17,20,22,25,28,35,36,39,42,43,46,48,49,50,\\
54,58,59,60,61,64,65,68,69,73,76,79,80,86,88,89,91,95,100,106\}.$
\\
\\
To prove that $R(5,14)\ge 239$, let $n=238$, and the set:\\
$S=\{3,8,9,11,12,13,15,17,20,21,25,27,32,36,37,42,45,49,52,54,\\
58,59,60,61,67,68,71,72,74,76,83,88,89,92,93,98,99,100,102,107,\\
108,119\}.$
\\
\\
To prove that $R(3,3,9)\ge 118$, let $n=117$, and the color sets:\\
Color $1: \{1,3,7,11,16,26,36,38,44,46,48,56\}$\\
Color $2: \{19,23,24,25,28,29,30,31,32,33,34,37,45\}$\\
Color $3: \{2,4,5,6,8,9,10,12,13,14,15,17,18,20,21,22,27,35,39,\\
40,41,42,43,47,49,50,51,52,53,54,55,57,58\}.$
\\
\\
To prove that $R(3,3,10)\ge 141$, let $n=140$, and the color sets:\\
Color $1: \{4,6,17,19,22,24,31,49,51,56,64,65,67\}$\\
Color $2: \{1,3,16,18,29,35,37,42,46,48,54,59,61,63,68\}$\\
Color $3: \{2,5,7,8,9,10,11,12,13,14,15,20,21,23,25,26,27,28,30,\\
32,33,34,36,38,39,40,41,43,44,45,47,50,52,53,55,57,58,60,62,66,\\
69,70\}.$
\\
\\
To prove that $R(3,3,11)\ge 158$, let $n=157$, and the color sets:\\
Color $1: \{3,4,16,22,24,30,36,45,51,57,62,63,68,74\}$\\
Color $2: \{6,7,9,10,23,26,28,31,39,42,50,53,58,61,66,69,77\}$\\
Color $3: \{1,2,5,8,11,12,13,14,15,17,18,19,20,21,25,27,29,32,33,\\
34,35,37,38,40,41,43,44,46,47,48,49,52,54,55,56,59,60,64,65,67,\\
70,71,72,73,75,76,78\}.$\\
\\
\\
By computer it is not hard to prove that the colorings are good.
\\
\\

\bf {References:}\\
\bf {[1]} Stanislaw P. Radziszowski: Small Ramsey Numbers (Dynamic Survey) 
\url{http://www.combinatorics.org/Surveys/ds1/sur.pdf}

\end{document}